\documentclass[review]{elsarticle}

\usepackage{hyperref}
\usepackage{natbib}
\usepackage{amsfonts,mathrsfs,upgreek,amsmath,amsthm}
\usepackage{bm}
\usepackage{color}

\journal{Journal of Statistical Planning and Inference}

\newtheorem{theorem}{Theorem}
\newtheorem{lemma}{Lemma}
\newtheorem{proposition}{Proposition}
\newtheorem{example}{Example}

\newtheorem{corollary}{Corollary}

\newcommand{\bpsi}{\mbox{\footnotesize\boldmath $\Uppsi$}}
\newcommand{\bbpsi}{\mbox{\small\boldmath $\Uppsi$}}

\biboptions{authoryear}

\begin{document}

\begin{frontmatter}

\title{On consistency of the MLE under finite mixtures of location-scale distributions with a structural parameter}

\author[add1]{Guanfu Liu} 

\author[add2]{Pengfei Li}

\author[add3]{Yukun Liu\corref{mycorrespondingauthor}}
\cortext[mycorrespondingauthor]{Corresponding author}
\ead{ykliu@sfs.ecnu.edu.cn}

\author[add3]{Xiaolong Pu}

\address[add1]{School of Statistics and Information, Shanghai University of International Business and Economics, Shanghai 201620, China}
\address[add2]{Department of Statistics and Actuarial Sciences, University of Waterloo, Waterloo, ON, Canada, N2L 3G1}
\address[add3]{School of Statistics,  East China Normal University,
                 Shanghai 200241, China}

\begin{abstract}
We provide a general and rigorous  proof for the strong  consistency of maximum likelihood estimators of
the cumulative distribution function of the mixing distribution and structural parameter
under finite mixtures of location-scale distributions with a structural parameter.
The consistency results do not require the parameter space of location and scale to be compact.
We illustrate the results by applying them to
finite  mixtures of location-scale distributions with the
component density function being one of the commonly used density functions: normal, logistic, extreme-value, or $t$.
An extension of the strong consistency results to finite mixtures of  multivariate elliptical distributions
is also discussed.

\end{abstract}


\end{frontmatter}


\section{Introduction}
Suppose we have an independent and identically distributed ($i.i.d.$) sample $X_1, \ldots, X_n$
from  the following finite mixture model:
\begin{equation}
\label{model}
g(x;\Psi,\sigma)= \sum_{j=1}^m\alpha_j f(x;\mu_j,\sigma)= \int_{\mathbb{R}} f(x;\mu,\sigma)d\Psi(\mu).
\end{equation}
Here $f(x;\mu,\sigma)$, the  component density function, is assumed to come from a location-scale distribution
family, namely, $f(x;\mu,\sigma)=\sigma^{-1}f\left((x-\mu)/\sigma;0,1\right)$ with $\mu \in \mathbb{R}$ and $\sigma \in  \mathbb{R}^+$
being the location and scale parameters, respectively.
The positive integer $m$
is called the order of the mixture model,
$(\alpha_1,\ldots,\alpha_m)$ with $\alpha_j\geq0$ and $\sum_{j=1}^m\alpha_j=1$
are called the mixing proportions,
and
$\Psi(\mu)=\sum_{j=1}^m\alpha_j I(\mu_j\leq \mu)$ is  called  the cumulative distribution function of the mixing distribution.
The parameter $\sigma$, appearing in all $m$ component density functions,
is called a structural parameter, and
Model (\ref{model}) is called a finite mixture of location-scale distributions
with a structural parameter.
Note that  $\Psi(\cdot)$ includes unknown $\mu_j$ and $\alpha_j$ parameters.
Hence, $(\Psi,\sigma)$ covers all the unknown parameters in (\ref{model}).
In this paper, we investigate the strong consistency of the maximum likelihood estimator (MLE) of $(\Psi,\sigma)$
under Model (\ref{model}).

Finite mixtures of
location-scale distributions with a structural parameter have many applications. They play
an important role in medical studies and genetics.
For example, \cite{Roeder:1994} applied the finite normal mixture model with a
structural parameter
to analyze sodium--lithium countertransport activity in red blood cells.
Finite mixtures of logistic distributions and of extreme value distributions
with  a structural parameter
are widely used to analyze failure-time data.
For instance,  a finite mixture of logistic distributions  with a structural parameter was used by \cite{Naya:2006} to study the thermogravimetric analysis trace.
A finite mixture of extreme value distributions
with a structural parameter was found to provide
an adequate fit to the logarithm of
the number of cycles to failure
for a group of 60 electrical appliances
(\citealp{Lawless:2003}; Example 4.4.2).
More applications can be found in  \cite{McLachlan:2000} and \cite{Lawless:2003}.

The maximum likelihood method  has been widely used to estimate
the unknown parameters in finite mixture models (\citealp{McLachlan:2000}; \citealp{Chen:2017}).
The consistency of the MLE under finite mixture models
has been studied by \cite{Kiefer:1956}, \cite{Redner:1981}, and \cite{Chen:2017}.
As pointed out by \cite{Chen:2017},
the results in \cite{Kiefer:1956}
require that {$g(x;\Psi,\sigma)$} can be continuously extended
to a compact space  of $(\Psi,\sigma)$.
This turns out to be impossible because $f(x;\mu,\sigma)$ is not well defined at $\sigma=0$.
To make the results in \cite{Kiefer:1956} applicable to our current setup,
we must constrain the parameter $\sigma$
to be in a compact subset of $\mathbb{R}^+$.
The consistency results in \cite{Redner:1981} require even more restrictive conditions:
the parameter space for $(\mu,\sigma)$
must be a compact subset of $\mathbb{R}\times \mathbb{R}^+$;
see \cite{Chen:2017} for more discussion.
It is worth mentioning that
\cite{Bryant:1991} established the strong consistency of
the estimators obtained by the linear-optimization-based method.
His result can be viewed as a generalization of the classical consistency result for MLE.
However, it requires that  the parameter space
be  closed  and  that  $m$  be equal to the true order of the mixture model.
By utilizing the properties of the normal distribution,
\cite{Chen:2017} proved the strong consistency of the MLE under
finite normal mixture models with a structural parameter without
imposing the compactness assumption on the parameter space.
To the best of our knowledge, general consistency results
for the MLE of $(\Psi,\sigma)$
under Model (\ref{model}) are not available in the literature except for the normal mixture model.

Because of the importance of finite mixtures of location-scale distributions with a structural parameter,
it is necessary   to study the consistency
of the MLE of the underlying parameters,  $(\Psi,\sigma)$, under Model (\ref{model}).
The goal of this paper is to   provide a general and rigorous proof of this consistency.
In Section 2, we present the main consistency results.
We emphasize that we do not require the parameter space of $(\mu,\sigma)$ to be compact.
The detailed proofs are given in Section 3. Section 4  illustrates  the consistency results by applying them
to  Model (\ref{model})  with $f(x;\mu,\sigma)$ being one of the commonly used component density
functions: normal, logistic, extreme-value, or $t$.
An extension of the consistency results to finite mixtures of multivariate elliptical distributions is discussed in Section 5.

\section{Main results}

With the $i.i.d.$ sample $X_1, \ldots, X_n$  from (\ref{model}),
the log-likelihood function of $(\Psi,\sigma)$ is given by
$$
\ell_n(\Psi,\sigma)=\sum_{i=1}^n \log\{ g(X_i;\Psi,\sigma)\}.
$$
The MLE of $(\Psi,\sigma)$ is defined as
$$
(\hat\Psi,\hat\sigma)=\arg\max_{\Psi\in {\bpsi}_m,~\sigma>0} \ell_n(\Psi,\sigma),
$$
where
$$
{\bbpsi}_m=\left\{\sum_{j=1}^m\alpha_j I(\mu_j\leq \mu):
\alpha_j\geq0,~\sum_{j=1}^m\alpha_j=1,~\mu_j\in \mathbb{R}
\right\}.
$$

In this section, we establish the consistency property of $(\hat \Psi,\hat\sigma)$ without imposing compactness on the parameter space of $(\mu,\sigma)$.
To discuss the consistency of $\hat\Psi$,
we define  \begin{equation}
\label{def.dist}
D\left(\Psi_1,\Psi_2\right)=\int_{\mathbb{R}}|\Psi_1(\mu)- \Psi_2(\mu) | \exp(-|\mu|)d\mu.
\end{equation}
We show that $D\left(\Psi_1,\Psi_2\right)$ is a distance on $\bbpsi_m$ in the Appendix.
Suppose $\Psi_0\in\bbpsi_m$ is the true cumulative distribution function of the mixing distribution.
We say that $\hat\Psi$ is strongly consistent if $D(\hat\Psi,\Psi_0)\to 0$ almost surely as $n\to\infty$.


The strong  consistency of  $(\hat \Psi,\hat\sigma)$ depends on
the following regularity conditions.

\begin{itemize}
\item[C1.] The finite mixture model in \eqref{model} is identifiable.
That is,  if $(\Psi_1,\sigma_1)$ and $(\Psi_2,\sigma_2)$
with $\Psi_1\in\bbpsi_m$, $\Psi_2\in\bbpsi_m$, $\sigma_1>0$, and $\sigma_2>0$
satisfy
$$
\int_{\mathbb{R}} f(x;\mu,\sigma_1)d\Psi_1(\mu)=\int_{\mathbb{R}} f(x;\mu,\sigma_2)d\Psi_2(\mu)
$$
for all $x$, then  $\Psi_1=\Psi_2$ and $\sigma_1=\sigma_2$.

\item[C2.] $\int_{\mathbb{R}} | \log  \{ g(x;\Psi_0,\sigma_0)\} |~  g(x;\Psi_0,\sigma_0) dx<\infty$,
where $(\Psi_0,\sigma_0)$ is the true value of $(\Psi,\sigma)$.

\item[C3.] There exist positive constants
$v_0$, $v_1$, and $\beta$ with $\beta>1$ such that for all $x$
$$
f(x;0,1)\leq \min \left\{v_0, v_1|x|^{-\beta}\right\}.
$$

\item[C4.]  The function  $f(x;0,1)$ is continuous with respect to $x$.
\end{itemize}

\begin{theorem}
\label{consistency}
Assume Conditions C1--C4 and that the true density is $g(x;\Psi_0,\sigma_0)$.
If an estimator $(\bar\Psi,\bar\sigma)$  of
$(\Psi,\sigma)$ satisfies
\begin{eqnarray}
\label{inequality}
\ell_n(\bar\Psi,\bar\sigma)-\ell_n(\Psi_0,\sigma_0)\geq c
>-\infty
\end{eqnarray}
for some constant $c$,
 then  $\bar\Psi \rightarrow \Psi_0$ with respect to the metric $D(\cdot,\cdot)$  in (\ref{def.dist})
 and $\bar\sigma\rightarrow\sigma_0$ almost surely as $n\rightarrow \infty$.
\end{theorem}

We comment that Conditions C1 and C2 are the standard regularity conditions
for the consistency of the MLE.
Condition C2 implies Wald's  integrability condition (\citealp{Wald:1949}).
Conditions C3 and C4 ensure that Lemma 1 in Section 3 is correct  and
the mixture density $g(x;\Psi,\sigma)$ can be continuously extended to
a compact space for $(\Psi,\sigma)$.
Conditions C2--C4 together guarantee that $\bar\sigma$ is away from 0 and  bounded above almost surely as $n\to\infty$.
Then the MLE consistency results in  \cite{Kiefer:1956} can be applied;
more discussion is given in Section 3.
In Section 4, we show that Model (\ref{model}) with $f(x;\mu,\sigma)$ being
one of
the four  commonly used  component density
functions (normal,  logistic, extreme-value, or $t$) satisfies these conditions.

It can be seen that the MLE $(\hat \Psi, \hat \sigma)$ satisfies (\ref{inequality}) with $c=0$.
Therefore,  by  Theorem  \ref{consistency},
both $\hat \Psi$ and $\hat \sigma$ are strongly consistent under Conditions C1--C4.

\begin{corollary}
\label{corollary}
Assume Conditions C1--C4 and that  the true density is $g(x;\Psi_0,\sigma_0)$.
Then  $\hat \Psi\rightarrow \Psi_0$ with respect to the metric $D(\cdot,\cdot)$ in (\ref{def.dist}) and $\hat\sigma\rightarrow \sigma_0$ almost surely as $n\rightarrow \infty$.
\end{corollary}



\section{Proofs}

\subsection{Some useful lemmas}
As discussed in \cite{Chen:2017},  except for Conditions C1 and C2,
a key regularity condition for  the MLE consistency results in  \cite{Kiefer:1956}
is that the definition of $g(x;\Psi,\sigma)$ can be
continuously extended  to  a compact space of $(\Psi,\sigma)$.
This extension could fail under our current setup
because $f(x;\mu,\sigma)$ is not well defined at $\sigma=0$.
To make the consistency results  in  \cite{Kiefer:1956} applicable,
a key step is to show that there exist positive constants
$\epsilon$ and $\Delta$ such that as $n\rightarrow\infty$,
the event sequence
$\{\epsilon\leq\bar\sigma\leq\Delta\}$ occurs almost surely.

We first present a technical lemma
that gives a uniform upper bound for the number of observations in the $\sigma^{1-a}$ neighborhood of $\mu$.
Here $a=(1+\beta)/(2\beta)$, where $\beta$ is given  in Condition C3.
With the condition that $\beta>1$ in Condition C3,  we have $0<a<1$.

For convenience of presentation, we let $b=2(\beta+1)/(\beta-1)$ and
$
\epsilon_0=\left(3mbv_0/\sigma_0\right)^{-1/(1-a)}.
$
Further, let
$$
G_n(x)=n^{-1}\sum_{i=1}^nI(X_i\leq x)
$$
be the empirical cumulative distribution function, and let $G(x)$ be the
cumulative distribution function of the $X_i$'s.
Define
\begin{equation*}
E=
\liminf\limits_{n\to\infty}\left\{\sup_{\mu\in\mathbb{R}} \sum_{i=1}^n I(|X_i-\mu|<\epsilon_0^{1-a})\leq n/(mb)\right\}.
\end{equation*}

%

\begin{lemma}
\label{lem1}
Suppose  $\{X_1,\ldots, X_n\}$  is an $i.i.d.$ sample from Model (\ref{model}).
Further, assume Conditions C3 and  C4 are satisfied.
Then we have
$$
P\left(E
\right )=1,
$$
and therefore almost surely  there exists an $n_0$
such that when $n\geq n_0$, we have
$$
  \sup_{\mu\in {\mathbb{R}}}  \sum_{i=1}^n I(|X_i-\mu|<\epsilon_0^{1-a})\leq   n/(mb).
$$
%
%
%
\end{lemma}

\proof

Since
$$
\sup_{\mu\in \mathbb{R}} \sum_{i=1}^n I(|X_i-\mu|<\epsilon_0^{1-a})
\leq \sup_{\mu\in \mathbb{R}}\hspace{0.01in} \big[ n \{ G_n(\mu+\epsilon_0^{1-a})-G_n(\mu-\epsilon_0^{1-a})\} \big],
$$
it suffices to show that
\begin{equation}
\label{lem1.eq1}
P\left(
\liminf\limits_{n\to\infty} \left\{ \sup_{\mu\in \mathbb{R}} \{
G_n(\mu+\epsilon_0^{1-a})-G_n(\mu-\epsilon_0^{1-a})
\}\leq 1/(mb)\right\}
\right )=1.
\end{equation}

Note that
\begin{eqnarray*}
\sup_{\mu\in\mathbb{R}}\left\{ G_n(\mu+\epsilon_0^{1-a})-G_n(\mu-\epsilon_0^{1-a}) \right\}
&\leq &\sup_{\mu\in\mathbb{R}}|  G_n(\mu+\epsilon_0^{1-a})-G(\mu+\epsilon_0^{1-a})|
\nonumber\\
&&+\sup_{\mu\in \mathbb{R}}|  G(\mu+\epsilon_0^{1-a})-G(\mu-\epsilon_0^{1-a})|\nonumber\\
&&+\sup_{\mu\in\mathbb{R}}|  G_n(\mu-\epsilon_0^{1-a})-G(\mu-\epsilon_0^{1-a})|.\nonumber
\label{empirical.dis.inequality}
\end{eqnarray*}
Under Conditions C3 and C4,
$$
\sup_{\mu\in\mathbb{R}}|  G(\mu+\epsilon_0^{1-a})-G(\mu-\epsilon_0^{1-a})|
\leq 2v_0\epsilon_0^{1-a}/\sigma_0.
$$

Let
$$
A_1=\left\{\lim_{n\to\infty} \sup_{\mu\in \mathbb{R}}|  G_n(\mu+\epsilon_0^{1-a})-G(\mu+\epsilon_0^{1-a})|=0\right\}
$$
and
$$
A_2=\left\{\lim_{n\to\infty} \sup_{\mu\in \mathbb{R}}|  G_n(\mu-\epsilon_0^{1-a})-G(\mu-\epsilon_0^{1-a})|=0\right\}.
$$
By the Glivenko--Cantelli theorem,
$$
P\left(A_1 \right)=P\left(A_2 \right)=1.
$$
Hence, $P(A_1\cap A_2)=1$.
Then  almost surely  there  exists an  $n_0$
such that when $n\geq n_0$,
\begin{eqnarray*}
\sup_{\mu\in \mathbb{R}}\left\{ G_n(\mu+\epsilon_0^{1-a})-G_n(\mu-\epsilon_0^{1-a}) \right\}
&\leq & 3 v_0\epsilon_0^{1-a}/\sigma_0=1/(mb),
\end{eqnarray*}
which implies (\ref{lem1.eq1}). This completes the proof.
\qed

The next lemma helps us to show that asymptotically we can restrict
$\sigma$ to be in a bounded interval away from 0 and $\infty$.
To  formally present the result, we define some notation.
Let $K_0=\int_{\mathbb{R}} \log \{ g(x;\Psi_0,\sigma_0)\} g(x;\Psi_0,\sigma_0)dx$. Condition C2 ensures that  $|K_0|<\infty$.
Further, define $\Delta=v_0/\exp(K_0-1)$.
We choose  a positive number $\epsilon$ that satisfies the following conditions:
\begin{itemize}
\item[D1.] $\epsilon\leq\epsilon_0=\left(3mbv_0/\sigma_0\right)^{-1/(1-a)}$;
\item[D2.] $\epsilon\leq (v_1/v_0)^{-2/(\beta+1)}$;
\item[D3.] $b^{-1}\log v_0+(1-b^{-1})\log v_1+\frac{1}{4}(\beta-1)\log\epsilon\leq K_0-1$.
\end{itemize}
Clearly, an $\epsilon$ satisfying the above conditions
 exists, since $\beta>1$ as assumed in Condition C3.

\begin{lemma}
\label{lemma2}
Suppose  $\{X_1,\ldots, X_n\}$  is   an $i.i.d.$ sample from Model (\ref{model}).
Further, assume Conditions  C1--C4 are satisfied.
Then we have
\begin{equation}
\label{lem2.ine1}
P\left(\lim_{n\to\infty}\left\{ \sup_{\Psi\in{\bpsi}_m, \sigma\in[\Delta,\infty)} \ell_n(\Psi,\sigma)
-\ell_n(\Psi_0,\sigma_0)\right\}=-\infty
\right )=1
\end{equation}
and
\begin{equation}
\label{lem2.ine2}
P\left(\lim_{n\to\infty}\left\{ \sup_{\Psi\in{\bpsi}_m, \sigma\in(0,\epsilon]} \ell_n(\Psi,\sigma)
-\ell_n(\Psi_0,\sigma_0)\right\}=-\infty
\right )=1.
\end{equation}

\end{lemma}

\proof

We start with (\ref{lem2.ine1}).
 Recall that  $f(x;\mu,\sigma)=\sigma^{-1}f((x-\mu)/\sigma; 0,1)$. By Condition C3, we have
\begin{eqnarray}
\ell_n(\Psi,\sigma)
&=&\sum_{i=1}^n\log\left\{\sum_{j=1}^m
\frac{\alpha_j}{\sigma} f\left(\frac{X_i-\mu_j}{\sigma};0,1\right)\right\}\nonumber\\
&\leq&\sum_{i=1}^n
\log\left\{\sum_{j=1}^m\frac{\alpha_j}{\sigma}v_0\right\}\nonumber\\
&\leq& n\left(\log v_0-\log \sigma\right),
\label{bound.of.likelihood}
\end{eqnarray}
where $v_0$ is given in Condition C3.

Hence, with $\Delta=v_0/\exp(K_0-1)$, we have
$$
\sup_{\Psi\in\bpsi_m,~\sigma\geq\Delta}\ell_n(\Psi,\sigma)-\ell_n(\Psi_0,\sigma_0)\leq n(K_0-1)-\ell_n(\Psi_0,\sigma_0).$$
By the strong law of large numbers and the definition of $K_0$, we have
(\ref{lem2.ine1}).

We next consider (\ref{lem2.ine2}).
Let $A=\{i: \min_{1\leq j\leq m} |X_{i}-\mu_j| \leq \sigma^{1-a}\}$  and $n(A)$ be the number of indices in set $A$.
For an index set  $S$, define
$
\ell_{n}(\Psi,\sigma;S)=\sum_{i\in S}\log g(X_i;\Psi,\sigma).
$
Similarly to (\ref{bound.of.likelihood}), it can be shown that
\begin{eqnarray}
\label{ln.lower1}
\ell_n(\Psi,\sigma;A)&\leq& n(A)(\log v_0-\log \sigma).
\end{eqnarray}
By Condition C3, we have
\begin{eqnarray*}
\ell_n(\Psi,\sigma;A^c)
&\leq&\sum_{i\in A^c}\log\left\{\sum_{j=1}^m
\frac{\alpha_j}{\sigma}v_1\left(\frac{|X_i-\mu_j|}{\sigma}\right)^{-\beta}\right\}\\
&=&\sum_{i\in A^c} \left[ \log v_1 + ( \beta-1) \log \sigma +
\log\left\{\sum_{j=1}^m
 \alpha_j   |X_i-\mu_j|^{-\beta}\right\} \right].
\end{eqnarray*}
Since
$ \min_{1\leq j\leq m} |X_{i}-\mu_j| \geq \sigma^{1-a}$  for all $i \in A^c$,
it follows that
\begin{eqnarray}
\ell_n(\Psi,\sigma;A^c)
  &\leq &\sum_{i\in A^c} \left\{  \log v_1 + ( \beta a-1) \log \sigma
   \right\}\nonumber \\
&=& n(A^c)\left\{ \log v_1+ \frac{1}{2} (\beta-1)\log \sigma \right\},
\label{ln.lower2}
\end{eqnarray}
where the   equality holds because  $\beta a=(\beta+1)/2>1$.

Combining (\ref{ln.lower1}) and (\ref{ln.lower2}) gives
\begin{align*}
~\ell_n(\Psi,\sigma)=&~\ell_n(\Psi,\sigma;A)+\ell_n(\Psi,\sigma;A^c)\nonumber\\
\leq&~  ~n(\log v_0-\log \sigma)+n(A^c)\left\{ \log v_1-\log v_0+  \frac{1}{2}(\beta+1) \log \sigma \right\}.\nonumber\\
=&~ n(\log v_0-\log \sigma)+n(1-b^{-1})\left\{ \log v_1-\log v_0+  \frac{1}{2}(\beta+1) \log \sigma \right\} \nonumber\\
&+\{n(A^c)-n(1-b^{-1})\} \left\{ \log v_1-\log v_0+  \frac{1}{2}(\beta+1) \log \sigma \right\}\nonumber \\
=&~n\left\{ b^{-1}\log v_0+(1-b^{-1})\log v_1+\frac{1}{4}(\beta-1)\log\sigma\right\}
\\
&+
\{n/b - n(A) \}
\left\{ \log v_1-\log v_0+  \frac{1}{2}(\beta+1) \log \sigma \right\},
\end{align*}
where in the last step we have used the fact that $b=2(\beta+1)/(\beta-1)$.
Then,   for $\sigma \in (0, \epsilon]$, since $\beta >1$ and $\epsilon$ satisfies D3, we have
\begin{align}
&\ell_n(\Psi,\sigma)-\ell_n(\Psi_0,\sigma_0)\nonumber\\
\nonumber&\leq n\left\{ b^{-1}\log v_0+(1-b^{-1})\log v_1+\frac{1}{4}(\beta-1)\log\epsilon\right\}-\ell_n(\Psi_0,\sigma_0)
\\
\nonumber&~~~+  \{n/b - n(A) \}  \left\{ \log v_1-\log v_0+  \frac{1}{2}(\beta+1) \log \sigma \right\}\\
&\leq n(K_0-1)-\ell_n(\Psi_0,\sigma_0)\label{ln.upper.part1}\\
&~~~+\{n/b - n(A) \}    \left\{ \log v_1-\log v_0+  \frac{1}{2}(\beta+1) \log \sigma \right\}
\label{ln.upper.part2}.
\end{align}
By the strong law of large numbers and the definition of $K_0$, we have
\begin{equation}
\label{ln.upper.part3}
P\left(
\lim_{n\to\infty}\left\{ n(K_0-1)-\ell_n(\Psi_0,\sigma_0)
\right\}=-\infty
\right)=1.
\end{equation}

Again for $\sigma\in(0, \epsilon]$,
since $\epsilon$ satisfies D2,   it follows that
\begin{equation}
\label{ln.upper.part4}
\log v_1-\log v_0+  \frac{1}{2}(\beta+1) \log \sigma
\leq
\log v_1-\log v_0+  \frac{1}{2}(\beta+1) \log \epsilon
\leq0.
\end{equation}
Since for each $i\in\{1,2,\dots,n\}$,
$$
I(\min_{1\leq j\leq m}|X_i-\mu_j|\leq\sigma^{1-a})\leq
\sum_{j=1}^m I(|X_i-\mu_j|\leq\sigma^{1-a}),
$$
we have
\begin{align*}
n(A)=&~\sum_{i=1}^nI(\min_{1\leq j\leq m}|X_i-\mu_j|\leq\sigma^{1-a})\\
\leq&~
\sum_{i=1}^n \sum_{j=1}^m I(|X_i-\mu_j|\leq\sigma^{1-a}) \\
\leq &~  m \cdot \sup_{\mu \in \mathbb{R}} \sum_{i=1}^n I(|X_i-\mu |\leq\sigma^{1-a}).
\end{align*}
Because $\sigma \in (0, \epsilon]$ and $\epsilon$ satisfies  D1,
it follows from Lemma 1 that
\begin{align*}
n(A)
\leq   m \cdot \sup_{\mu \in \mathbb{R}} \sum_{i=1}^n I(|X_i-\mu |\leq\epsilon_0^{1-a})
\leq  n/b.
\end{align*}
This together with \eqref{ln.upper.part4} implies that  for large enough $n$,
$$
\{n/b-n(A)\} \left\{ \log v_1-\log v_0+  \frac{1}{2}(\beta+1) \log \sigma \right\}
\leq 0,
$$
almost surely.
Combining (\ref{ln.upper.part1})--(\ref{ln.upper.part4}) leads to (\ref{lem2.ine2}). This completes the proof.
\qed

\subsection{Proof of Theorem \ref{consistency}}

%

The results in (\ref{lem2.ine1}) and (\ref{lem2.ine2})
imply that
$$
P\Big(\liminf\limits_{n\to\infty}\{\epsilon\leq \bar\sigma\leq \Delta\} \Big)=1.
$$
Hence, we can confine $\bar\sigma$  to $[\epsilon,\Delta]$ asymptotically,
and the results of \cite{Kiefer:1956} apply.
For completeness, we outline the key steps  of the proof of Theorem 1 by following \cite{Kiefer:1956}.

In the first step, we compactify the parameter space $\bbpsi_m$.
Let  $$\bar \bbpsi_m=\{\gamma+\rho \Psi: \Psi\in \bbpsi_m,~\gamma\geq 0, ~\rho\geq 0, ~0\leq\gamma+\rho\leq1\}.$$
We extend the definition of $D(\cdot,\cdot)$ to $\bar \bbpsi_m$ without modifications.
Then $\bar \bbpsi_m$ is a compact metric space with respect to $D(\cdot,\cdot)$.
See  the Appendix for a proof of its compactness.
Let $\bm{\theta}=(\Psi,\sigma)$
and define
$$
\bar D(\bm\theta_1,\bm\theta_2)
=D(\Psi_1,\Psi_2)+|\sigma_1-\sigma_2|.
$$
Since $D(\cdot,\cdot)$ is a distance
and
$$
D(\Psi_1,\Psi_2)\leq 2\int_{\mathbb{R}}\exp(-|\mu|)d\mu=4,
$$
we can verify that $\bar D(\cdot,\cdot)$ is a bounded distance
on $\bar \bbpsi_m\times [\epsilon,\Delta]$.
Hence, $\bar \bbpsi_m\times [\epsilon,\Delta]$ is compact with respect to $\bar D(\cdot,\cdot)$.

In the second step, we argue that $g(x;\Psi,\sigma)$
is continuous for all $x$ on $\bar \bbpsi_m\times [\epsilon,\Delta]$ under the distance
$\bar D(\cdot,\cdot)$.
That is,  for any $\bm\theta=(\Psi,\sigma)\in \bar \bbpsi_m\times [\epsilon,\Delta]$,
if $\bar D(\tilde{\bm\theta},\bm\theta)\to0$, then
we have
$$
g(x;\tilde\Psi,\tilde\sigma)\to g(x;\Psi,\sigma).
$$
For the given $(\Psi,\sigma)$, we define
$$
H(\mu,\sigma^*)=\Psi(\mu) I(\sigma\leq\sigma^*)
$$
and define $\tilde H$  to be $H$ with  $(\tilde \Psi, \tilde \sigma)$ in place of $(  \Psi,   \sigma)$.
Then the mixture density  can be written as
$$
g(x;\Psi,\sigma)=  \int_{\mathbb{R}\times\mathbb{R}^+} f(x;\mu,\sigma^*)dH(\mu,\sigma^*).
$$
We further define
$$
\tilde D(H_1,H_2)
=\int_{\mathbb{R}\times\mathbb{R}^+} |H_1(\mu,\sigma^*)-H_2(\mu,\sigma^*) | \exp\{-|\mu|-|\sigma^*|\}d\mu d\sigma^*.
$$
Lemma 2.4 of
\cite{Chen:2017} implies that
if $f(x;\mu,\sigma)$ is continuous for $(\mu,\sigma)$,
$\lim_{|\mu|+|\sigma|\to\infty} f(x;\mu,\sigma)=0$,
and $\tilde D(\tilde H,H)\to0$, then
we have $g(x;\tilde\Psi,\tilde\sigma)\to g(x;\Psi,\sigma)$.

Condition C4 ensures that  $f(x;\mu,\sigma)$ is continuous for $(\mu,\sigma)$;
Condition C3 ensures that $\lim_{|\mu|+|\sigma|\to\infty} f(x;\mu,\sigma)=0$.
Hence, to argue that $g(x;\Psi,\sigma)$
is continuous for all $x$, we need to show only that if  $\bar D(\tilde{\bm\theta},\bm\theta)\to0$, then we must have
$\tilde D(\tilde H,H)\to0$. Note that
\begin{eqnarray}
\nonumber\tilde D(\tilde H,H)
&=&\int_{\mathbb{R}\times\mathbb{R}^+} | \tilde\Psi(\mu)I(\tilde\sigma\leq\sigma^*)- \Psi(\mu)I( \sigma\leq\sigma^*) | \exp\{-|\mu|-|\sigma^*|\}d\mu d\sigma^*\\
\nonumber&\leq& \int_{\mathbb{R}\times\mathbb{R}^+} | \tilde\Psi(\mu) - \Psi(\mu) | \exp\{-|\mu|-|\sigma^*|\}d\mu d\sigma^*\\
\nonumber&&+\int_{\mathbb{R}\times\mathbb{R}^+} | I(\tilde\sigma\leq\sigma^*)- I( \sigma\leq\sigma^*) | \exp\{-|\mu|-|\sigma^*|\}d\mu d\sigma^*\\
\nonumber&\leq &\int_{\mathbb{R}} | \tilde\Psi(\mu) - \Psi(\mu) | \exp\{-|\mu|\}d\mu \\
\nonumber&&+2\int_{\mathbb{R}^+} | I(\tilde\sigma\leq\sigma^*)- I( \sigma\leq\sigma^*) | \exp\{- |\sigma^*|\}d\sigma^*\\
&\leq&D(\tilde\Psi,\Psi)+2|\tilde\sigma-\sigma|\exp\{-\min(\tilde\sigma,\sigma)\}.
\label{tilde.dist}
\end{eqnarray}
Further, when $\bar D(\tilde{\bm\theta},\bm\theta)\to0$, we must have $D(\tilde\Psi,\Psi)\to0$ and $\tilde\sigma\to\sigma$,
which, together with (\ref{tilde.dist}), implies that $\tilde D(\tilde H,H)\to0$.
Hence, $g(x;\Psi,\sigma)$
is continuous for all $x$ on $\bar \bbpsi_m\times [\epsilon,\Delta]$ under the distance
$\bar D(\cdot,\cdot)$.

In the third step, we show that for any $\bm\theta=(\Psi,\sigma)\in \bar \bbpsi_m\times [\epsilon,\Delta]$ such that $(\Psi,\sigma)\neq (\Psi_0,\sigma_0)$,
we can find a $\delta$ such that
\begin{equation}
\label{step3.obj}
E\left\{
\log g\big(X;B_{\delta} (\bm\theta)\big)
\right\}
-{E\left\{\log
g\big(X;\Psi_0,\sigma_0\big)
\right\}}<0,
\end{equation}
where
$$
B_{\delta} (\bm\theta)=\Big\{\bm\theta^*=(\Psi^*,\sigma^*):  \bar D(\bm\theta^*,\bm\theta )<\delta,~~\bm\theta^*\in \bar \bbpsi_m\times [\epsilon,\Delta]\Big\}
$$
and
$$
g\big(x; B_{\delta} (\bm\theta)\big)=\sup_{ {\bm\theta^*}\in B_{\delta} (\bm\theta) } g(x;\Psi^*,\sigma^*).
$$
By  \cite{Kiefer:1956} and \cite{Chen:2017} and because of Conditions C1, C2, and C4,
we can obtain {(\ref{step3.obj})} by arguing that
\begin{equation}
\label{step3.upper}
\lim_{\delta\to0^+}
E\left[
\left\{
\log\frac{g\big(X;B_{\delta} (\bm\theta)\big) }{g\big(X;\Psi_0,\sigma_0\big)
}
\right\}^+
\right]<\infty.
\end{equation}
Here $x^+=\max(x,0)$.
By Condition C3,
$$
\log\frac{g\big(X;B_{\delta} (\bm\theta)\big) }{g\big(X;\Psi_0,\sigma_0\big)
}\leq \log (v_0/\epsilon_0)-\log g\big(X;\Psi_0,\sigma_0\big).
$$
Hence,
$$
\left\{
\log\frac{g\big(X;B_{\delta} (\bm\theta)\big) }{g\big(X;\Psi_0,\sigma_0\big)
}
\right\}^+\leq |\log (v_0/\epsilon_0)|+|\log g\big(X;\Psi_0,\sigma_0\big)|,
$$
which, together with Condition C2, implies (\ref{step3.upper}). Hence, (\ref{step3.obj}) is proved.

In the last step, we combine the above results to complete the proof of Theorem 1.
For any $\delta>0$, let
$$
B^{(0)}=B_{\delta}(\bm\theta_0).
$$
The result in  (\ref{step3.obj}) leads to a finite open cover of
the compact set $\{\bar\bbpsi_m\times [\epsilon,\Delta]\}\backslash B^{(0)}$.
Then by the finite covering theorem,  we can find
 $\bm\theta_1,\ldots,\bm\theta_K$ from $\bar\bbpsi_m\times [\epsilon,\Delta]$  and positive $\delta_1,\ldots,\delta_K$
for some positive integer $K$ such that
$$
\{\bar \bbpsi_m\times [\epsilon,\Delta]\}\backslash B^{(0)}\subset \bigcup_{k=1}^K B_{\delta_k}(\bm\theta_k)
$$
and
\begin{equation}
\label{final.ine}
E\left\{
\log g\big(X;B_{\delta_k} (\bm\theta_k)\big)
\right\}
-E\left\{\log
g\big(X;\Psi_0,\sigma_0\big)
\right\}<0.
\end{equation}
Define
$$
B^{(k)}=B_{\delta_k} (\bm\theta_k),~~k=1,\ldots,K,~~
B^{(K+1)}=\{ \bar\bbpsi_m\times (0,\epsilon] \}\bigcup \{ \bar\bbpsi_m\times  [\Delta,\infty)  \}.
$$
Then the whole parameter space  of $\bm\theta=(\Psi,\sigma)$ can be written as
$$
\bar\bbpsi_m\times(0,\infty)=\bigcup_{k=0}^{K+1} B^{(k)}.
$$

Note that it can be easily verified that
$$
 \sup_{\Psi\in{\bpsi}_m, \sigma\in(0,\epsilon]} \ell_n(\Psi,\sigma)
 = \sup_{\Psi\in{\bar \bpsi}_m, \sigma\in(0,\epsilon]} \ell_n(\Psi,\sigma)
$$
and
$$
 \sup_{\Psi\in{\bpsi}_m, \sigma\in[\Delta,\infty)} \ell_n(\Psi,\sigma)
 = \sup_{\Psi\in{\bar \bpsi}_m, \sigma\in[\Delta,\infty)} \ell_n(\Psi,\sigma).
$$
Hence, by the strong law of large numbers, (\ref{lem2.ine1})--(\ref{lem2.ine2}), and (\ref{final.ine}), we have for $k=1,\ldots,K+1$
\begin{equation*}
P\left(\lim_{n\to\infty}\left\{ \sup_{{\bm\theta}\in B^{(k)}}  \ell_n(\Psi,\sigma)
-\ell_n(\Psi_0,\sigma_0)\right\}=-\infty
\right )=1.
\end{equation*}
Therefore, for any estimator $(\bar\Psi,\bar\sigma)$ of
$(\Psi,\sigma)$ such that
\begin{eqnarray*}
\ell_n(\bar\Psi,\bar\sigma)-\ell_n(\Psi_0,\sigma_0)\geq c
>-\infty
\end{eqnarray*}
for some constant $c$, we must have for any $\delta>0$
$$
P\Big(\liminf\limits_{n\to\infty}\{\bar{\bm\theta}\in B^{(0)} \} \Big)=
P\Big(\liminf\limits_{n\to\infty}\{\bar D(\bar{\bm\theta}, {\bm\theta_0})<\delta\} \Big)=1.
$$
This implies that $\bar D(\bar{\bm\theta}, {\bm\theta_0})\to0$ almost surely.
Equivalently, $D(\bar\Psi,\Psi_0)\to0$ and $\bar\sigma\to\sigma_0$ almost surely.
This completes the proof of Theorem 1.

\section{Examples}
In this section,
 we illustrate the consistency results of Section 2
by showing  that the four  commonly used   component density
functions satisfy  Conditions C1--C4.
Consequently, by Corollary \ref{corollary}, the MLE of $(\Psi,\sigma)$ is strongly consistent if
$f(x; \mu, \sigma)$ is one of these  four functions.
As preparation, we present a sufficient condition for verifying Condition C2.
\begin{proposition}
\label{prop1}
If $\int_{\mathbb{R}}  \{\log f(x;\mu,1)\}  f(x;0,1) dx>-\infty$ for any given $\mu$ and
Condition C3 is satisfied,
then Condition C2 is satisfied.
\end{proposition}
\proof
Let
$$
\Psi_0(\mu)=\sum_{j=1}^{m_0}\alpha_{0j}I(\mu_{0j}\leq\mu).
$$
Without loss of generality, we assume $\sigma_0=1$.
Otherwise, we can take the transformation $t=x/\sigma_0$.

Note that
\begin{eqnarray*}
&&\int_{\mathbb{R}} | \log  \{ g(x;\Psi_0,\sigma_0)\} |~  g(x;\Psi_0,\sigma_0) dx\\
&=&
\int_{\mathbb{R}} \left|
\log\left\{
\sum_{j=1}^{m_0} \alpha_{0j} f(x;\mu_{0j},1)
\right\}
\right|~\sum_{j=1}^{m_0} \alpha_{0j} f(x;\mu_{0j},1) dx\\
&=&
\int_{\mathbb{R}} \left|
\log\left\{
\sum_{j=1}^{m_0} \alpha_{0j} \frac{f(x;\mu_{0j},1)}{v_0}
\right\}+\log v_0
\right|~\sum_{j=1}^{m_0} \alpha_{0j} f(x;\mu_{0j},1) dx\\
&\leq& |\log v_0|+ \int_{\mathbb{R}} \left|
\log\left\{
\sum_{j=1}^{m_0} \alpha_{0j} \frac{f(x;\mu_{0j},1)}{v_0}
\right\}
\right|~\sum_{j=1}^{m_0} \alpha_{0j} f(x;\mu_{0j},1) dx\\
&=&|\log v_0|-\int_{\mathbb{R}}
\log\left\{
\sum_{j=1}^{m_0} \alpha_{0j} \frac{f(x;\mu_{0j},1)}{v_0}
\right\}
~\sum_{j=1}^{m_0} \alpha_{0j} f(x;\mu_{0j},1) dx,
\end{eqnarray*}
where in the last step we have used the fact that $f(x;0,1)\leq v_0$ in Condition C3.
Hence, to verify Condition C2, it suffices to show that
\begin{equation}
\label{equi.condi.c2}
\int_{\mathbb{R}}
\log\left\{
\sum_{j=1}^{m_0} \alpha_{0j} f(x;\mu_{0j},1)
\right\}
~\sum_{j=1}^{m_0} \alpha_{0j} f(x;\mu_{0j},1) dx>-\infty.
\end{equation}
By Jensen's inequality,
\begin{eqnarray*}
&&\int_{\mathbb{R}}
\log\left\{
\sum_{j=1}^{m_0} \alpha_{0j} f(x;\mu_{0j},1)
\right\}
~\sum_{j=1}^{m_0} \alpha_{0j} f(x;\mu_{0j},1) dx\\
&\geq&
\int_{\mathbb{R}}
\left\{
\sum_{j=1}^{m_0} \alpha_{0j}
\log
 f(x;\mu_{0j},1)
\right\}
~\sum_{j=1}^{m_0} \alpha_{0j} f(x;\mu_{0j},1) dx\\
&=&\sum_{j=1}^{m_0} \sum_{h=1}^{m_0}  \alpha_{0j}\alpha_{0h}
\int_{\mathbb{R}}\{\log
 f(x;\mu_{0j},1)\} f(x;\mu_{0h},1) dx\\
 &=&\sum_{j=1}^{m_0} \sum_{h=1}^{m_0}  \alpha_{0j}\alpha_{0h}
\int_{\mathbb{R}}\{\log
 f(t;\mu_{0j}-\mu_{0h},1)\} f(t;0,1) dt.
\end{eqnarray*}
By the condition   $\int_{\mathbb{R}}  \{ \log f(x;\mu,1) \}  f(x;0,1) dx>-\infty$ for any given $\mu$,
we have (\ref{equi.condi.c2}).
This completes the proof.
\qed

Note that it is easy to verify that the four commonly used component density functions all satisfy  Condition C4.
We  now verify that they all satisfy Conditions C1--C3.

\begin{example}[Normal distribution]
Let $f(x;0,1)=(2\pi)^{-1/2}\exp(-x^2/2)$, the probability density function of the standard normal distribution.
\end{example}

\begin{itemize}
\item[(a)] The identifiability of Model \eqref{model} with the normal  component density function
 follows from \cite{Teicher:1963} and \cite{Yakowitz:1968}.  Hence,
 Condition C1 is satisfied.
\item[(b)] Condition C3 is satisfied if  we choose $v_0=v_1=\left( 2\pi\right)^{-1/2}$ and $\beta=2$.
\item[(c)] For any given $\mu$,
$$
\int_{\mathbb{R}}  \{\log f(x;\mu,1)\}  f(x;0,1) dx
={-0.5\log(2\pi)-0.5(\mu^2+1)>-\infty.}
$$
Thus,
Condition C2 is also satisfied by Proposition \ref{prop1}.
\end{itemize}

\vspace{.2cm}
\begin{example}[Logistic distribution]
Let $f(x;0, 1)=e^x/(1+e^x)^2$, the probability density function  of the standard logistic distribution.
\end{example}

\begin{itemize}
\item[(a)] Following Theorem 2.1 of \cite{Holzmann:2004}, Model \eqref{model}
with a logistic  component density function is identifiable. Hence,
 Condition C1 is satisfied.

\item[(b)] It can be verified that Condition C3 is satisfied with $v_0=v_1=1$ and $\beta=2$.
\item[(c)]
Since $\log (1+e^x)\leq \log 2+|x|\leq \log 2+0.5(1+x^2)$ for any $x$, it follows that
\begin{eqnarray*}
\int_{\mathbb{R}}  \{\log f(x;\mu,1)\}  f(x;0,1) dx
&=&-\mu-2\int_{\mathbb{R}}\left\{\log\left(1+e^{x-\mu}\right)\right\} f(x;0,1) dx\\
&\geq&-\mu-2\log 2-\int_{\mathbb{R}}\{1+(x-\mu)^2\} f(x;0,1) dx\\
&=&-\mu-2\log 2-(1+\mu^2+\pi^2/3)>-\infty.
\end{eqnarray*}
Thus,
Condition C2 is also satisfied by Proposition \ref{prop1}.
\end{itemize}

\vspace{.2cm}
\begin{example}[Extreme-value distribution]
Let $f(x;0,1)=\exp\{-x-\exp(-x)\}$,  the probability density function
of the standard extreme value type I distribution or the Gumbel distribution.
\end{example}
\begin{itemize}
\item[(a)] The identifiability of Model \eqref{model} with the  component density function being the probability density function of the Gumbel distribution
follows from \cite{Ahmad:2010}. Hence,
 Condition C1 is satisfied.

\item[(b)] It can be verified that Condition C3 is satisfied with $v_0=v_1=1$ and $\beta=2$.

\item[(c)]
Note that
\begin{eqnarray*}
\int_{\mathbb{R}}  \{\log f(x;\mu,1)\}  f(x;0,1) dx
&=& \int_{\mathbb{R}} \{-(x-\mu)-e^{-(x-\mu)}\} \exp(-x-e^{-x})dx\\
&=&\mu-\gamma-e^{\mu}>-\infty.
\end{eqnarray*}
Here $\gamma$ is the Euler--Mascheroni constant.
Hence, Condition C2 is satisfied by Proposition \ref{prop1}.

\end{itemize}

\vspace{.2cm}
\begin{example}[Student's $t$ distribution]
Let $f(x; 0, 1)= (1+x^2/\nu)^{-(\nu+1)/2} C_{\nu}$,  the probability density function of Student's $t$ distribution
with $\nu$ degrees of freedom. Here $\nu$ is a given positive integer,
$C_{\nu} =\Gamma\{(\nu+1)/2\} \{ \sqrt{\nu \pi} \Gamma (\nu/2) \}^{-1} $,
and $\Gamma(\cdot)$ is the Gamma function.
\end{example}

\begin{itemize}
\item[(a)] The identifiability of Model \eqref{model} with the  component density function being the probability density function of the $t$-distribution
follows from \cite{Holzmann:2006}. Hence,
 Condition C1 is satisfied.

\item[(b)] It can be verified that Condition C3 is satisfied with $v_0=1, v_1=\nu$, and $\beta=2$.
\item[(c)] Note that
\begin{eqnarray*}
&&\int_{\mathbb{R}}  \{\log f(x;\mu,1)\}  f(x;0,1) dx\\
&=&\log C_{\nu}-\frac{\nu+1}{2} C_{\nu} \int_{\mathbb{R}} \log \left\{1+\frac{(x-\mu)^2}{\nu}\right\}  (1+x^2/\nu)^{-(\nu+1)/2}  dx.\end{eqnarray*}
By the comparison test for improper integrals,
$$
0\leq \int_{\mathbb{R}} \log \left\{1+\frac{(x-\mu)^2}{\nu}\right\}  (1+x^2/\nu)^{-(\nu+1)/2}  dx<\infty
$$ for any given $\mu$.
Hence, $\int_{\mathbb{R}}  \log f(x;\mu,1)  f(x;0,1) dx>-\infty$
and Condition C2 is satisfied by Proposition \ref{prop1}.

\end{itemize}

\section{Extension to multivariate case}

In this section, we extend the results in Corollary  \ref{corollary} to
finite mixtures of  multivariate elliptical distributions, a special class of finite mixtures of multivariate location-scale distributions.
The identifiability of this special class of models has been well studied in  \cite{Holzmann:2006}.

Suppose we have  $i.i.d.$  $p$-dimensional random vectors $\bm X_1, \ldots, \bm X_n$
from  the following finite mixture model:
\begin{equation}
\label{model2}
g(\bm x;\Psi,\bm\Sigma)= \sum_{j=1}^m\alpha_j f(\bm x;\bm\mu_j,\bm\Sigma)= \int_{\mathbb{R}^p} f(\bm x;\bm\mu, \bm\Sigma)d\Psi(\bm\mu).
\end{equation}
Here $f(\bm x;\bm\mu, \bm\Sigma)$, the  component density function, is assumed to
take the form
$$
f(\bm x;\bm\mu, \bm\Sigma)=|\bm\Sigma|^{-1/2}f_0\left((\bm x-\bm\mu)^\tau\bm\Sigma^{-1}(\bm x-\bm\mu)\right),$$
where $\bm x$, $\bm\mu \in \mathbb{R}^p$,  $\bm\Sigma$ is a $p\times p$ positive definite matrix, and $f_0(x)$ is a density generator,
i.e., a non-negative function on $[0,\infty)$  such that $f_0(\bm x^\tau\bm x)$ is a probability density function.
The  MLE $(\hat\Psi,\bm{\hat\Sigma})$ of $(\Psi, \bm{\Sigma})$  is defined as for the univariate case.
Next we extend the distance $D(\cdot,\cdot)$ from the univariate to the multivariate case:
\begin{equation}
\label{def.dist2}
D^*(\Psi_1,\Psi_2)=\int_{\mathbb{R}}|\Psi_1(\bm\mu)- \Psi_2(\bm\mu)| \exp(-|\bm\mu|)
d{\bm \mu}.
\end{equation}
Here for $\bm\mu=(\mu_1,\ldots,\mu_p)^\tau$, $|\bm\mu|$ is interpreted as $\sum_{l=1}^p|\mu_l |$.
Similarly to the proof for $D(\cdot,\cdot)$ in the Appendix,
we can verify that $D^*(\cdot,\cdot)$ is a distance.

The consistency of $(\hat\Psi,\bm{\hat\Sigma})$ relies on the following regularity conditions:
\begin{itemize}
\item[C1*.] The finite mixture model in \eqref{model2} is identifiable.

\item[C2*.] $\int_{\mathbb{R}^p} | \log  \{ g(\bm x;\Psi_0,\bm\Sigma_0)\} |   g(\bm x;\Psi_0,\bm\Sigma_0)  d{\bm x}<\infty$,
where $(\Psi_0,\bm\Sigma_0)$ is the true value of $(\Psi,\bm\Sigma)$.
\item[C3*.] There exist positive constants
$v_0$, $v_1$, and $\beta$ with $\beta>p$ such that for all $x\geq 0$
$$
f_0(x)\leq \min \left\{v_0, v_1 x^{-\beta/2}\right\}.
$$

\item[C4*.]  For $x\geq 0$, $f_0(x)$ is continuous in $x$.
\end{itemize}

Under the regularity conditions C1*--C4*, we have the strong consistency of $(\hat\Psi,\bm{\hat\Sigma})$
in the following theorem.

\begin{theorem}
\label{thm2}
Assume Conditions C1*--C4* and that the true density is $g(x;\Psi_0,\bm\Sigma_0)$.
Then  $\hat \Psi\rightarrow \Psi_0$  with respect to the metric $D^*(\cdot,\cdot)$ in (\ref{def.dist2})
and $\hat{\bm\Sigma}\rightarrow \Sigma_0$ almost surely as $n\rightarrow \infty$.
\end{theorem}

One of the key steps  of the proof of Theorem \ref{thm2} is to establish a similar result to Lemma \ref{lem1}
for the multivariate case.
This result can be obtained by combining the proof of Corollary 3 in \cite{Alexandrovich:2014} with that of Lemma \ref{lem1}.
See the Appendix for a proof of Theorem \ref{thm2}.

As an illustration, we consider finite mixtures of multivariate normal distributions with a common and unknown variance-covariance matrix $\bm\Sigma$.
\begin{example} (Multivariate normal distribution).
Let $f_0(x)=(2\pi)^{-1/2}\exp(-x/2)$, the density generator for the multivariate normal density function.
Clearly, $f_0(x)$  satisfies Condition C4*.

\begin{itemize}
\item[(a)] The identifiability of finite mixtures of multivariate normal distributions is covered by \cite{Holzmann:2006}.  Hence,
 Condition C1* is satisfied.

\item[(b)] It can be verified that Condition C3* is satisfied with $v_0=(2\pi)^{-p/2}$,
$$v_1=(2\pi)^{-p/2}(p+1)^{(p+1)/2},$$  and $\beta=p+1$.
\item[(c)]
Following the proof of Proposition \ref{prop1}, to verify Condition C2*, it suffices
to show that
$$
\int_{\mathbb{R}^p}\{\log  f(\bm x;\bm\mu, \bm\Sigma_0)\} f(\bm x; \bm 0,\Sigma_0)d\bm x>-\infty
$$
for all $\bm\mu$.  Note that
\begin{eqnarray*}
&& \int_{\mathbb{R}^p}\{\log  f(\bm x;\bm\mu, \bm\Sigma_0)\} f(\bm x; \bm 0,\Sigma_0)d\bm x \\
&=&-0.5p\log(2\pi)-0.5\log|\bm\Sigma_0|
-2p-\bm\mu^\tau{\bm\Sigma_0}^{-1}\bm\mu  \\
&>&-\infty.
\end{eqnarray*}
Hence, Condition C2* is satisfied.
\end{itemize}
\end{example}

Since finite mixtures of multivariate normal distributions with a common and unknown variance-covariance matrix satisfy regularity conditions C1*--C4*,
the MLE $(\hat\Psi,\bm{\hat\Sigma})$ under this model is strongly consistent.

\section{Summary and discussion}
In this paper, we establish the strong consistency of the cumulative distribution function of
the mixing distribution
and the structural parameter in finite mixtures of location-scale distributions with a structural parameter
in both the univariate and multivariate cases.
We further demonstrate that some commonly used finite mixtures of location-scale distributions
satisfy the regularity conditions.

For the model setups in (\ref{model}) and (\ref{model2}), $\Psi$
is assumed to have finite support, and
the scale parameter $\sigma$ or $\bm\Sigma$ is assumed to be the same in all the component density functions.
Two considerations  underlie these assumptions.
First,  if $\Psi$  is fully nonparametric, then the mixture model may not be identifiable.
For example, the normal mixture model is not identifiable if $\Psi$ is fully nonparametric \citep{Chen:2017}.
Assuming that $\Psi$ has finite support ensures that
finite mixtures of some commonly used location-scale distributions are
identifiable.
Second, if the $\sigma$ or $\bm\Sigma$ can vary in different component density functions,
then the log-likelihood is unbounded \citep{Chen:2008,Chen:2009}. Hence, the usual MLEs of the unknown parameters are not well defined.
We may need to consider other estimation methods such as the penalized MLE considered in \cite{Chen:2008}, \cite{Chen:2009}, and \cite{Alexandrovich:2014}.
We leave the consistency properties of such  estimators
to future research.

We next discuss the applicability of Corollary \ref{corollary} and Theorem \ref{thm2}.
The results in Corollary \ref{corollary} and Theorem \ref{thm2} are applicable only to the MLE, i.e., the global maximum point of the log-likelihood function. Commonly used algorithms such as the EM-algorithm may lead to
a local maximum point of the log-likelihood, which is not guaranteed  by our result to be consistent.
In practice, we suggest trying multiple initial values to increase the chance of locating the global maximum point.

\section*{A. Appendix}
\subsection*{{\rm A.1.} Proof that $D\left(\cdot,\cdot\right)$ is a distance}

To show that $D\left(\cdot,\cdot\right)$ is a distance.,  it suffices to show
\begin{itemize}
\item[(a)] $D\left(\Psi_1,\Psi_2\right)\geq 0$;
\item[(b)] $D\left(\Psi_1,\Psi_3\right)\leq D\left(\Psi_1,\Psi_2\right)+D\left(\Psi_2,\Psi_3\right)$ for any $\Psi_1,\Psi_2, \Psi_3\in \bbpsi_m$;
\item[(c)] $D\left(\Psi_1,\Psi_2\right)=D\left(\Psi_2,\Psi_1\right)$;
\item[(d)] $D\left(\Psi_1,\Psi_2\right)=0$ if and only if $\Psi_1(\mu)=\Psi_2(\mu)$ for all $\mu\in \mathbb{R}$.
\end{itemize}

Based on the definition of $D(\cdot,\cdot)$ in (\ref{def.dist}), it is easy to verify that (a)--(c) are satisfied.
Next we discuss (d).
If $\Psi_1(\mu)=\Psi_2(\mu)$ for all $\mu\in \mathbb{R}$, then obviously
$$
D\left(\Psi_1,\Psi_2\right)=\int_{\mathbb{R}}|\Psi_1(\mu)- \Psi_2(\mu)| \exp(-|\mu|)d\mu=0.
$$
We now argue that if  $D\left(\Psi_1,\Psi_2\right)=0$, we must have $\Psi_1(\mu)=\Psi_2(\mu)$ for all $\mu\in \mathbb{R}$. We denote the distinct values of the supports of $\Psi_1(\mu)$ and $\Psi_2(\mu)$ as $\{t_1<t_2<\cdots<t_{m^*}\}$
and define $t_{m^*+1}=+\infty$.
 We write $\Psi_1(\mu)$ and $\Psi_2(\mu)$ as
$$
\Psi_1(\mu)=\sum_{j=1}^{m^*}\alpha_{j1}I(t_j\leq\mu)
\mbox{
and
}
\Psi_2(\mu)=\sum_{j=1}^{m^*}\alpha_{j2}I(t_j\leq\mu).
$$
If $\Psi_1(\mu)\neq \Psi_2(\mu)$ for some $\mu$, then we can find a $j_0$ such that
$$
\alpha_{j1}=\alpha_{j2}\mbox{ for }j\leq j_0-1~~\mbox{ and}~~\alpha_{j_01}\neq\alpha_{j_02}.
$$
This implies that
\begin{eqnarray*}
D\left(\Psi_1,\Psi_2\right)
&\geq &\int_{\mu\in [t_{j_0},t_{j_0+1}) } |\Psi_1(\mu)-\Psi_2(\mu)| \exp(-|\mu|)d\mu\\
&\geq &|\alpha_{j_01}-\alpha_{j_02}| \int_{t_{j_0}}^{t_{j_0+1} } \exp(-|\mu|)d\mu>0.
\end{eqnarray*}
Hence, if  $D\left(\Psi_1,\Psi_2\right)=0$, we must have $\Psi_1(\mu)=\Psi_2(\mu)$ for all $\mu\in \mathbb{R}$.
This completes the proof.

\subsection*{{\rm A.2.} Proof of the compactness of $\bar\bbpsi_m$}

We prove that $\bar\bbpsi_m$ is compact with respect to the distance $D(\cdot,\cdot)$ according to
the following  equivalent definition of compactness of a metric space.
A metric space
is compact if and only if every sequence in this space has a convergent subsequence
whose limit is also in this space.

In the proof of the compactness of $\bar\bbpsi_m$, we need the following
results from  real analysis.
\begin{enumerate}
\item[] {\it Result (i)}. If $\{c_n,~n=1,2,\dots\}$ is a bounded real sequence, then $\{c_n,~n=1,2,\dots\}$
has a convergent subsequence $\{c_{n_k},~k=1,2,\dots\}$, and its limit is finite.

\item[] {\it Result (ii)}.  If $\{c_n,~n=1,2,\dots\}$ is a unbounded real sequence, then $\{c_n,~n=1,2,\dots\}$
has a subsequence $\{c_{n_k},~k=1,2,\dots\}$ that diverges to $\infty$ or $-\infty$.
\end{enumerate}

Result (i) is just the classic Bolzano--Weierstrass theorem from real analysis.
Hence, we  give only the proof for Result (ii).
Suppose $\{c_n,~n=1,2,\dots\}$ has no upper bound.
Then for any given positive integer $k$,
there exists $n_k$ such that $c_{n_k}>k$.
Then $\lim_{k\to \infty}c_{n_k}=\infty$. Hence, Result (ii) holds.

Next we return to the proof of the compactness of $\bar\bbpsi_m$.
Let $\{\Psi_{L},L=1,2,\dots\}$
be a sequence in $\bar\bbpsi_m$,
where
$$
\Psi_{L}(\mu)=\gamma_L+\rho_L\sum_{j=1}^m\alpha_{jL}I(\mu_{jL}\leq \mu)
$$
with
$\gamma_L\geq 0$, $\rho_L\geq 0$, $0\leq\gamma_L+\rho_L\leq 1$,
$\alpha_{jL}\geq 0$, $\sum_{j=1}^m\alpha_{jL}=1$, and
$-\infty<\mu_{1L}\leq\mu_{2L}\leq\ldots\leq\mu_{mL}<\infty $.

Using Results (i)--(ii) and G. Cantor's ``diagonal method,"
we can find a subsequence $\{L_k,k=1,2,\ldots\}$
such that
$$\gamma_0=\lim_{k\to\infty}\gamma_{L_k}, ~
\rho_0=\lim_{k\to\infty}\rho_{L_k},~
\alpha_{j0}=\lim_{k\to\infty}\alpha_{jL_k},~
\mu_{j0}=\lim_{k\to\infty}\mu_{jL_k},
$$
where $\gamma_0\geq 0$, $\rho_0\geq 0$, $0\leq\gamma_0+\rho_0\leq 1$,
$\alpha_{j0}\geq 0$, $\sum_{j=1}^m\alpha_{j0}=1$, and
$-\infty\leq\mu_{10}\leq\mu_{20}\leq\ldots\leq\mu_{m0}\leq\infty $.
We define
$$\Psi_{0}(\mu)=\gamma_0+\rho_0\sum_{j=1}^m\alpha_{j0}I(\mu_{j0}\leq \mu).$$
Further, we define two index sets $S_1$ and $S_2$ as
$S_1=\{j:~\mu_{j0}=-\infty\}$ and $S_2=\{j:~\mu_{j0}=\infty\}$, respectively. Let $S=S_1\cup S_2$.
Then $\Psi_{0}(\mu)$ can be rewritten as
$$
\Psi_{0}(\mu)=\gamma_0'+\rho'_0\sum_{j\notin S}\alpha_{j0}'I(\mu_{j0}\leq \mu),
$$
where $\gamma_0'=\gamma_0+\rho_0\sum_{j\in S_1}\alpha_{j0}$,
$\rho'_0=\rho_0\sum_{j\notin S}\alpha_{j0}$, and
$\alpha_{j0}'=\alpha_{j0}/\sum_{l\notin S}\alpha_{l0}$ for $j\notin S$.
This implies that $\Psi_{0}\in\bar\bbpsi_m$.

Let $\Psi_{L_k}=\gamma_{L_k}+\rho_{L_k}\sum_{j=1}^m\alpha_{jL_k}I(\mu_{jL_k}\leq \mu)$.
To finish the proof, we need to argue that as $k\to\infty$,
$$
D(\Psi_{L_k},\Psi_{0})\to 0.
$$
Note that
\begin{align}
\nonumber D(\Psi_{L_k},\Psi_0)=&~\int_{\mathbb{R}}|\Psi_{L_k}(\mu)-\Psi_0(\mu)|\exp(-|\mu|)d\mu\\
\nonumber=&\int_{\mathbb{R}}|\gamma_{L_k}-\gamma_{0}+\sum_{j=1}^m\left\{\rho_{L_k}\alpha_{jL_k}I(\mu_{jL_k}\leq \mu)
-\rho_{0}\alpha_{j0}I(\mu_{j0}\leq \mu)\right\}|
\exp(-|\mu|)d\mu\\
\leq& \int_{\mathbb{R}}\left|\gamma_{L_k}-\gamma_{0}\right|\exp(-|\mu|)d\mu\label{limit.d1}\\
&+\sum_{j=1}^m\int_{\mathbb{R}}|\rho_{L_k}\alpha_{jL_k}I(\mu_{jL_k}\leq \mu)
-\rho_{0}\alpha_{j0}I(\mu_{j0}\leq \mu)|
\exp(-|\mu|)d\mu\label{limit.d2}.
\end{align}

Since $\gamma_0=\lim_{k\to\infty}\gamma_{L_k}$,
we have that as $k\to\infty$,
\begin{equation}
\label{limit.d11}
\int_{\mathbb{R}}\left|\gamma_{L_k}-\gamma_{0}\right|\exp(-|\mu|)d\mu
=2\left|\gamma_{L_k}-\gamma_{0}\right|\to 0.
\end{equation}
By the triangular inequality and the facts that $0\leq \rho_0\leq 1$ and $0\leq \alpha_{j0}\leq 1$,
we have
\begin{align}
&\int_{\mathbb{R}}|\rho_{L_k}\alpha_{jL_k}I(\mu_{jL_k}\leq \mu)
-\rho_{0}\alpha_{j0}I(\mu_{j0}\leq \mu)|
\exp(-|\mu|)d\mu\nonumber\\
&\leq\int_{\mathbb{R}}|\rho_{L_k}\alpha_{jL_k}
-\rho_{0}\alpha_{j0}|
\exp(-|\mu|)d\mu\label{limit.d21}\\
&+\int_{\mathbb{R}}|
I(\mu_{jL_k}\leq \mu)-I(\mu_{j0}\leq \mu)|
\exp(-|\mu|)d\mu\label{limit.d22}
\end{align}
Similarly to (\ref{limit.d11}),  as $k\to\infty$ we have
\begin{equation}
\label{limit.d111}
\int_{\mathbb{R}}|\rho_{L_k}\alpha_{jL_k}
-\rho_{0}\alpha_{j0}|
\exp(-|\mu|)d\mu \to 0.
\end{equation}
For (\ref{limit.d22}),
we have that if $\mu_{j0}$ is finite,
\begin{align}
\int_{\mathbb{R}}|
I(\mu_{jL_k}\leq \mu)-I(\mu_{j0}\leq \mu)|
\exp(-|\mu|)d\mu\leq |\mu_{jL_k}- \mu_{j0}|\to 0\label{limit.d221}
\end{align}
as $k\to\infty$, and   if $\mu_{j0}=\infty$ or $-\infty$,
\begin{align}
\int_{\mathbb{R}}|
I(\mu_{jL_k}\leq \mu)-I(\mu_{j0}\leq \mu)|
\exp(-|\mu|)d\mu =\exp(-|\mu_{jL_k}|)\to 0
\label{limit.d222}
\end{align}
as $k\to\infty$.

Combining (\ref{limit.d21})--(\ref{limit.d222}) leads to
$$
\lim _{k\to\infty} \int_{\mathbb{R}}|\rho_{L_k}\alpha_{jL_k}I(\mu_{jL_k}\leq \mu)
-\rho_{0}\alpha_{j0}I(\mu_{j0}\leq \mu)|
\exp(-|\mu|)d\mu=0,
$$
which, together with (\ref{limit.d1})--(\ref{limit.d11}), implies that
$$
\lim_{k\to\infty} D(\Psi_{L_k},\Psi_0)=0.
$$
This completes the proof.

\subsection*{{\rm A.3.} Proof of Theorem \ref{thm2}}

Our proof of Theorem \ref{thm2} is similar to that of Theorem \ref{consistency}.
Hence, we  simply  outline the key steps.

In the first step, we establish a result
similar to Lemma \ref{lem1}.
Let $a^*=(\beta+p)/(2\beta)$ and  $b^*=2(\beta-p+2)/(\beta-p)$, where $\beta>p$ is given in Condition C3*.
Hence, $0<a^*<1$ and $b^*>1$.
We choose $\epsilon_0^*$ such that
$$
\frac{v_0\pi^{p/2}}{|\bm\Sigma_0|^{1/2}\Gamma(p/2+1)}(\epsilon_0^*)^{ (1-a^*)/2}= 1/( 2mb^*).
$$
For the matrix $\bm\Sigma$,  ``$\bm\Sigma>0$" means that $\bm\Sigma$ is a positive definite matrix.

\begin{lemma}
\label{lemS1}
Suppose  $\{\bm X_1,\ldots, \bm X_n\}$  is an $i.i.d.$ sample from $g(\bm x;\Psi_0,\bm\Sigma_0)$.
Further, assume Condition C3*.
Let
$$
E^*=
\liminf\limits_{n\to\infty}\left\{\sup_{\{\bm\Sigma>0, |\bm\Sigma|\leq\epsilon_0^*\}}\sup_{\bm\mu\in \mathbb{R}^p}
\sum_{i=1}^n I \left\{(\bm x_i-\bm\mu)^\tau\bm\Sigma^{-1}(\bm x_i-\bm\mu)\leq {|\bm\Sigma|}^{-a^*/p}\right\}\leq n/(mb^*)\right\}.
$$
Then $P(E^*)=1$.
\end{lemma}
\proof
Following the   proof of Corollary 3 in \cite{Alexandrovich:2014},
we can show that
$$
P\left(
\liminf\limits_{n\to\infty}\left\{\sup_{\{\bm\Sigma>0, |\bm\Sigma|\leq\epsilon_0^*\}}\sup_{\bm\mu\in \mathbb{R}^p}
\sum_{i=1}^n I \left\{(\bm x_i-\bm\mu)^\tau\bm\Sigma^{-1}(\bm x_i-\bm\mu)\leq {|\bm\Sigma|}^{-a^*/p}\right\}\leq a_n
\right\} \right)=1,
$$
where
$$
a_n=\frac{3}{4}\sqrt{n\log\log n}+\frac{nv_0\pi^{p/2}}{|\bm\Sigma_0|^{1/2}\Gamma(p/2+1)}(\epsilon_0^*)^{(1-a^*)/2}
=\frac{3}{4}\sqrt{n\log\log n}+ \frac{n}{ 2mb^*}.
$$
When $n$ is large enough, $a_n\leq n/(mb^*)$.
Hence,
$
P\left(
E^* \right)=1.
$
This completes the proof.
\qed

In the second step, we establish a result similar to Lemma \ref{lemma2}.
We first define some notation.
Based on  $i.i.d.$  $p$-dimensional random vectors ${\bm X_1}, \ldots, {\bm X_n}$ from $g(\bm x;\Psi,\bm\Sigma)$,
the log-likelihood of $(\Psi,\bm\Sigma)$ is
$$
\ell_n(\Psi,\sigma)=\sum_{i=1}^n \log\{ g(X_i;\Psi,\bm\Sigma)\}.
$$
Let $K_0^*=\int_{\mathbb{R}^p} \log \{ g(\bm x;\Psi_0,\bm\Sigma_0)\} g(\bm x;\Psi_0,\bm\Sigma_0)d\bm x$. Condition C2* ensures that  $|K_0^*|<\infty$.
Further, define $\Delta^*=\{ v_0/\exp( K_0^*-1)\}^2$.
We choose a positive number $\epsilon^*$ that satisfies the following conditions:
\begin{itemize}
\item[D1*.] $\epsilon^*\leq\epsilon_0^*$;
\item[D2*.] $\epsilon^*\leq (v_1/v_0)^{-4/(\beta-p+2)}$;
\item[D3*.] $(1/b^*)\log v_0+(1-1/b^*)\log v_1+\frac{1}{8}(\beta-p)\log\epsilon^*\leq K_0^*-1$.
\end{itemize}
Clearly, an $\epsilon^*$ satisfying the above conditions  exists, since $\beta>p$ as assumed in Condition C3*.

\begin{lemma}
\label{lemmaS2}
Suppose  $\{X_1,\ldots, X_n\}$  is an $i.i.d.$ sample from $g(\bm x;\Psi_0,\bm\Sigma_0)$.
Further, assume Conditions C1*--C4* are satisfied.
Then we have
\begin{equation}
\label{lemS2.ine1}
P\left(\lim_{n\to\infty}\left\{ \sup_{\Psi\in{\bpsi}_m, \bm\Sigma>0,|\bm\Sigma|\geq\Delta^*} \ell_n(\Psi,\bm\Sigma)
-\ell_n(\Psi_0,\bm\Sigma_0)\right\}=-\infty
\right )=1
\end{equation}
and
\begin{equation}
\label{lemS2.ine2}
P\left(\lim_{n\to\infty}\left\{ \sup_{\Psi\in{\bpsi}_m, \bm\Sigma>0,|\bm\Sigma|\leq\epsilon^*} \ell_n(\Psi,\bm\Sigma)
-\ell_n(\Psi_0,\bm\Sigma_0)\right\}=-\infty
\right )=1.
\end{equation}

\end{lemma}

\proof

We start with (\ref{lemS2.ine1}).
By Condition C3*, we have
\begin{eqnarray}
\ell_n(\Psi,\bm\Sigma)
&=&\sum_{i=1}^n\log\left\{\sum_{j=1}^m
\frac{\alpha_j}{|\bm\Sigma|^{1/2}} f_0\left((\bm X_i-\bm\mu_j)^\tau\bm\Sigma^{-1}(\bm X_i-\bm\mu_j)\right)\right\}\nonumber\\
&\leq& n\left(\log v_0-0.5\log |\bm\Sigma|\right),
\label{app.bound.of.likelihood}
\end{eqnarray}
where $v_0$ is given in Condition C3*.

Hence, with $\Delta^*=\{ v_0/\exp( K_0^*-1)\}^2$, we have
$$
\sup_{\Psi\in{\bpsi}_m, \bm\Sigma>0,|\bm\Sigma|\geq\Delta^*} \ell_n(\Psi,\bm\Sigma)
-\ell_n(\Psi_0,\bm\Sigma_0)\leq n(K_0^*-1)-\ell_n(\Psi_0,\bm\Sigma_0).$$
By the strong law of large numbers and the definition of $K_0^*$, we have
(\ref{lemS2.ine1}).

We next consider (\ref{lemS2.ine2}).
Let
$$
A^*=\left\{i: \min_{1\leq j\leq m} (\bm X_i-\bm\mu_j)^\tau\bm\Sigma^{-1}(\bm X_i-\bm\mu_j)\leq {|\bm\Sigma|}^{-a^*/p} \right\}.
$$
Similarly to (\ref{app.bound.of.likelihood}), it can be shown that
\begin{eqnarray}
\label{app.ln.lower1}
\ell_n(\Psi,\bm\Sigma;A^*)&\leq& n(A^*)(\log v_0-0.5\log| \bm\Sigma|).
\end{eqnarray}
By Condition C3*, we have
\begin{eqnarray}
\nonumber
\ell_n(\Psi,\bm\Sigma;(A^*)^c)
&\leq&\sum_{i\in (A^*)^c}\log\left[\sum_{j=1}^m
\frac{\alpha_j}{|\bm\Sigma|^{1/2}}v_1\left\{(\bm X_i-\bm\mu_j)^\tau\bm\Sigma^{-1}(\bm X_i-\bm\mu_j) \right\}^{-\beta/2}\right]\\
\nonumber&\leq&\ n\big((A^*)^c\big)\left\{ \log v_1+0.5(a^*\beta/p-1)\log |\bm\Sigma|\right\}\\
&=& n\big((A^*)^c\big)\left\{ \log v_1+0.25(\beta-p)\log |\bm\Sigma|\right\},
\label{app.ln.lower2}
\end{eqnarray}
where in the last step we have used the fact that $a^*\beta=(\beta+p)/2$.

Combining (\ref{app.ln.lower1}) and (\ref{app.ln.lower2}) gives
\begin{align*}
~\ell_n(\Psi,\bm\Sigma)=&~\ell_n(\Psi,\bm\Sigma;A^*)+\ell_n(\Psi,\bm\Sigma;(A^*)^c)\nonumber\\
\leq&~n(\log v_0-0.5\log |\bm\Sigma|)+n\big((A^*)^c\big) \left\{ \log v_1-\log v_0+  0.25(\beta-p+2) \log |\bm\Sigma| \right\}\nonumber\\
=&~ n(\log v_0-0.5\log |\bm\Sigma|)+n(1-1/b^{*})\left\{ \log v_1-\log v_0+  0.25(\beta-p+2) \log |\bm\Sigma| \right\} \nonumber\\
&+\left\{n\big((A^*)^c\big)-n(1-1/b^{*})\right\} \left\{ \log v_1-\log v_0+  0.25(\beta-p+2) \log |\bm\Sigma| \right\}\nonumber \\
=&~n\left\{ (1/b^*)\log v_0+(1-1/b^*)\log v_1+\frac{1}{8}(\beta-p)\log|\bm\Sigma|\right\}
\\
&+ \{n/b^* - n\big( A^* \big)  \} \left\{ \log v_1-\log v_0+  0.25(\beta-p+2) \log |\bm\Sigma|  \right\},
\end{align*}
where in the last step we have used the fact that $b^*=2(\beta-p+2)/(\beta-p)$.
Then,  for $\bm\Sigma$ satisfying $|\bm\Sigma|\leq\epsilon^*$,
since $\epsilon^*$ satisfies D3*,  we have
\begin{align}
&\ell_n(\Psi,\bm\Sigma)-\ell_n(\Psi_0,\bm\Sigma_0)\nonumber\\
&\leq n(K_0^*-1)-\ell_n(\Psi_0,\bm\Sigma_0)\label{app.ln.upper.part1}\\
&~~~+ \{n/b^* - n\big( A^* \big)  \} \left\{ \log v_1-\log v_0+  0.25(\beta-p+2) \log |\bm\Sigma|  \right\}
\label{app.ln.upper.part2}.
\end{align}

By the strong law of large numbers and the definition of $K_0^*$, we have
\begin{equation}
\label{app.ln.upper.part3}
P\left(
\lim_{n\to\infty}\left\{ n(K_0^*-1)-\ell_n(\Psi_0,\bm\Sigma_0)
\right\}=-\infty
\right)=1.
\end{equation}
Again  for $\bm\Sigma$ satisfying $|\bm\Sigma|\leq\epsilon^*$,
since $\epsilon^*$ satisfies D2*, we have
\begin{equation}
\label{app.ln.upper.part4}
\log v_1-\log v_0+  0.25(\beta-p+2) \log |\bm\Sigma|
\leq
\log v_1-\log v_0+  0.25(\beta-p+2) \log \epsilon^*
\leq 0,
\end{equation}
which together with Lemma \ref{lemS1} implies that, for large enough $n$,
$$
\{n/b^* - n\big( A^* \big)\}  \left\{ \log v_1-\log v_0+  0.25(\beta-p+2) \log |\bm\Sigma|  \right\}
\leq 0,
$$
almost surely.
Combining (\ref{app.ln.upper.part1})--(\ref{app.ln.upper.part4}) leads to (\ref{lemS2.ine2}). This completes the proof.
\qed

Lemma \ref{lemmaS2} implies that asymptotically we can confine $(\hat\Psi,\hat{\bm\Sigma})$
to $\bm\Theta$, where $\bm\Theta=\bbpsi_m\times\{\bm\Sigma:\bm\Sigma>0,~|\bm\Sigma|\in[\epsilon^*,\Delta^*]\}$.
Note that $\bm\Theta$ is completely regular \citep{Chen:2009}, so
the techniques in \cite{Wald:1949} and \cite{Kiefer:1956} can be applied to establish
the strong consistency of $(\hat\Psi,\hat{\bm\Sigma})$.
This completes the proof of Theorem 2.

\section*{Acknowledgements}
We are grateful to the editor, the associate editor,  and three anonymous referees for their
insightful and constructive comments which led to an improved presentation of this article.
The work of Dr.~Liu and Professor Pu is supported by grants from the National Natural Science Foundation
of China (11771144, 11371142, 11771145,  11471119), the
outstanding doctoral dissertation cultivation plan of action (PY2015049),  the Program of Shanghai Subject Chief Scientist
(14XD1401600), and the 111 Project (B14019).
Dr.~Li's research is supported in part by NSERC Grant RGPIN-2015-06592.

\section*{Reference}
\bibliography{mybibfile}

\end{document}